\documentclass[12pt, letter]{amsart}

\usepackage{geometry}

\usepackage[colorlinks = true,
            linkcolor = red,
            urlcolor  = magenta,
            citecolor = black,
            anchorcolor = blue]{hyperref}
\usepackage{color}
\usepackage{upgreek}
\usepackage{tikz}
\usepackage[applemac]{inputenc} % for Macs
\usetikzlibrary{arrows}
\usepackage{xcolor}
\usepackage[all]{xy}
\usepackage{graphicx}
\usepackage{wrapfig}
\usepackage{yfonts}
\usepackage{graphicx}
\usepackage{amssymb}
\usepackage{epstopdf}
\usepackage{mathrsfs}
\usepackage[mathcal]{eucal}
\usepackage{bm}

\newtheorem{Thm}{Theorem}[section]

\newtheorem{Prop}[Thm]{Proposition}
\newtheorem{Cor}[Thm]{Corollary}

\theoremstyle{definition}
\newtheorem{Def}[Thm]{Definition}
\theoremstyle{remark}

\newtheoremstyle{named}{}{}{\itshape}{}{\bfseries}{.}{.5em}{#1 #3}
\theoremstyle{named}

\def\si{{\Sigma}}

\def\Z{\mathbb{Z}}

\def\g{\mathfrak{g}}
\def\Frenkel:2013uda{\mathfrak{h}}

\def\cO{\mathcal{O}}

\def\cQ{\mathcal{Q}}

\def\c{\phi}

\def\d{\delta}

\def\t{\tau}

\def\w{\omega}

\def\bo{\textbf{o}}

\def\=>{\Longrightarrow}

\def\to{\longrightarrow}

\def\o+{\oplus}
\def\bo+{\bigoplus}

\def\<{\langle}
\def\>{\rangle}
\def\({\left(}
\def\){\right)}

\def\^{\wedge}
\def\+{\dagger}

\def\inv{^{-1}}

\def\dd[#1,#2]{\frac{d#1}{d#2}}
\def\del[#1,#2]{\frac{\partial #1}{\partial #2}}
\def\over[#1]{\overline{#1}}

\def\tab{\;\;\;\;\;\;}

\makeatletter
\def\mr@ignsp#1 {\ifx\:#1\@empty\else #1\expandafter\mr@ignsp\fi}%
\newcommand{\multiref}[1]{\begingroup%\let\protect\string%
\xdef\mr@no@sparg{\expandafter\mr@ignsp#1 \: }%
\def\mr@comma{}%
\@for\mr@refs:=\mr@no@sparg\do{\mr@comma\def\mr@comma{,}\ref{\mr@refs}}%
\endgroup}
\makeatother

\newcommand{\hypref}[2]{\ifx\href\asklFrenkel:2013udaas #2\else\href{#1}{#2}\fi}

\newcommand{\case}[2][cccccccccccccccccccccccccccccccccccccccccc]{\left\{\begin{array}{#1}#2 \\ \end{array}\right.}
\newcommand{\Eq}[1]{\begin{align}#1\end{align}}

\tikzset{->-/.style={decoration={
  markings,
  mark=at position .5 with {\arrow{latex}}},postaction={decorate}}}
\tikzset{
    >=latex
    }

\newcommand{\nc}{\newcommand}
\nc{\on}{\operatorname}
\nc{\la}{\lambda}
\nc{\wh}{\widehat}
\nc{\ghat}{\wh\g}
\nc{\mb}{\mathbf}

\begin{document}

\title{Super Riemann Surfaces and Fatgraphs}

\author[A.S.  Schwarz]{Albert S. Schwarz}
\address[Albert S. Schwarz]
{Department of Mathematics,
University of California at Davis,
Davis, CA, USA\newline
Email: \href{mailto:schwarz@math.ucdavis.edu}{schwarz@math.ucdavis.edu}
}

\author[A.M. Zeitlin]{Anton M. Zeitlin}
\address[Anton M. Zeitlin]{
          Department of Mathematics, 
          Louisiana State University, 
          Baton Rouge, LA, USA\newline
Email: \href{mailto:zeitlin@lsu.edu}{zeitlin@lsu.edu},\newline
 \href{http://math.lsu.edu/~zeitlin}{http://math.lsu.edu/$\sim$zeitlin}
}

\date{\today}

\begin{abstract}
Our goal is to describe superconformal structures on super Riemann 
surfaces (SRS), based on data assigned to a fatgraph. 
We start from the complex structures on punctured $(1|1)$-supermanifolds, characterizing the corresponding moduli and the deformations using Strebel differentials and certain \v{C}ech cocycles for a specific covering, which we reproduce from a fatgraph data, consisting of $U(1)$-graph connection and odd parameters at the vertices. Then we consider dual $(1|1)$-supermanifolds and related superconformal structures for $N=2$ super Riemann surfaces. The superconformal structures $N=1$ SRS are computed as the fixed points of involution on supermoduli space of $N=2$ SRS.
\end{abstract}

\numberwithin{equation}{section}

\maketitle

\setcounter{tocdepth}{1}
%\tableofcontents
\numberwithin{equation}{section}

\addtocontents{toc}{\protect\setcounter{tocdepth}{1}}
\addtocontents{toc}{\protect\setcounter{tocdepth}{2}}

\section{Introduction}\label{sec:intro}
\subsection{Some history and earlier results}
The geometry of moduli spaces of (punctured) Riemann surfaces has been a central topic in modern mathematics for many years. Since the 1980s, string theory served as a significant source of ideas in studying moduli spaces. For a proper description of string theory, one has to consider certain generalizations of moduli spaces related to the fact that strings, while propagating, should carry extra anticommutative parameters, thus generating what is known as  superconformal manifold as introduced by  M.A. Baranov and A.S. Schwarz \cite{baranov1985multiloop} or super Riemann surface (SRS) as independently introduced by D. Friedan \cite{friedan1986notes} (see also \cite{baranov1989geometry}, \cite{manin1986critical}, \cite{giddings1988geometry}, \cite{CR},  \cite{schwarz}, and \cite{superwitten} for a review). It turned out that such spaces' geometry is quite involved, e.g., \cite{donagi}. 
An important task is, of course, related to the parametrization of such supermoduli.

There are several ways of looking at the parametrization problem. For example, one could deal with supermoduli spaces of punctured Riemann surfaces with the negative Euler characteristic from the point of view of higher Teichm\"uller theory as a subset in the character variety for the corresponding supergroup. In the case of original moduli spaces using methods of hyperbolic geometry R. Penner described coordinates in the universal cover of moduli space, the Teichm\"uller space, as the subspace of the character variety of $PSL(2,\mathbb{R})$,
so that the corresponding Riemann surfaces appear here from the uniformization point of view as a factor of the upper half-plane by the element of the related character variety, i.e., the Fuchsian subgroup \cite{pbook}.

The action of the mapping class group in these coordinates is rational. It could be described combinatorially using decorated triangulations or dual objects, known as metric fatgraphs or ribbon graphs for the corresponding Riemann surfaces. 
Thus constructed coordinates were generalized to the case of reductive groups \cite{fg}.  The supergroup case yet remained a mystery until recently. In \cite{penzeit}, \cite{IPZ}, \cite{IPZ2}, such coordinates were constructed in the framework of the higher Teichm\"uller spaces associated to supergroups $OSp(1|2)$ and $OSp(2|2)$, which correspond to the Teichm\"uller spaces $N=1$ and $N=2$ SRS. The desired $N=1$ and $N=2$ SRS could be reconstructed using the elements of character variety via the appropriately modified uniformization approach \cite{CR}, \cite{natanzon2004moduli}.

There is a different, more ``hands-on" approach to the moduli spaces of punctured Riemann surfaces, where one can see directly the transition functions for the corresponding complex structures, which we discuss in more detail below. One can start from the parameterization of moduli spaces by the so-called Strebel differentials, which again can be described using metric fatgraphs \cite{konts}. That approach allowed (see \cite{mulase}) to ``glue" the Riemann surface explicitly by constructing transition functions. 

In this paper, we want to generalize this construction in the case of super Riemann surfaces. We start by describing the moduli space of  $(1|1)$-supermanifolds. This result also describes the moduli space of
  $N=2$ super Riemann surfaces. Finally, we study the moduli space of  $N=1$ super Riemann surfaces using the fact that this space can be obtained as a set of fixed points of the involution  of the space of $(1|1)$-supermanifolds
  constructed in \cite{schwarz}.

We would also like to mention recent progress in studying supermoduli spaces from various perspectives. While our approach deals with real parametrization of supermoduli in parallel to work on super-Teichm\"uller theory \cite{penzeit},\cite{IPZ},\cite{IPZ2}, a lot of exciting features of supermoduli are related to the algebro-geometric description. The main result of Donagi and Witten \cite{DW} that supermoduli space is not projected led to a renewed interest in the subject in the modern era and supergeometry in general. One can mention recent works of Felder, Kazhdan, and Polishchuk dealing with the Schottky approach for supermoduli \cite{superperiod} as well as the general treatment of supermoduli spaces as Deligne-Mumford stacks \cite{DM}. Some other recent results, which use both real and complex geometry points of view, are related to enumerative invariants related to supermoduli \cite{PN}, \cite{Yau}.

\subsection{The structure of the paper and main results.}
In Section 2 we review basic notions related to $(1|1)$-supermanifolds, $N=1$ and $N=2$ Super-Riemann surfaces (SRS). We devote special attention to the punctured $N=1$ SRS with two puncture classes corresponding to various spin structure choices: Ramond (R) and Neveu-Schwarz (NS).

In Section 3, we define two instrumental objects which come from geometric topology.
The first object is a fatgraph (or ribbon graph). This graph is homotopically equivalent to the punctured surface with the cyclic ordering of half-edges at every vertex, 
which comes from the orientation of the surface so that each puncture is associated with 
a particular cycle on the graph. 
The second object is a spin structure on the fatgraph, making it a {\it spin fatgraph}. We describe spin structures as the classes of orientations on fatgraphs based on the works \cite{penzeit}, \cite{IPZ}, \cite{IPZ2}, where $N=1$, $N=2$ SRS were studied from a uniformization perspective. This construction allows distinguishing boundary components of such spin fatgraph, separating them into two sets based on comparing their orientation and the orientation induced by the surface. Those two sets correspond to NS and R punctures in the uniformization picture. 

Section 4 is devoted to an important construction allowing us to relate the data assigned to the fatgraphs to the theory of moduli of Riemann Surfaces, following \cite{mulase}, \cite{konts}.  
 Namely, we explicitly describe the moduli spaces of moduli spaces of surfaces 
 $F^c$ with marked points, using special covering $\{U_v, V_p\}$, with one neighborhood $U_v$ for every vertex $v$ and $V_p$ for every puncture $p$. The set of $\{U_v\}$ has only double overlaps $U_v\cap U_{v'}$, corresponding to edges $\{v,v'\}$, so that $\cup_{v}U_v=F$ is a punctured surface. 
 $U_p$ overlaps with all $U_v$ for all the vertices surrounding the puncture. 
  To construct the corresponding transition functions $w'=f_{v',v}(w), y=f_{p,v}(w)$ on overlaps,  we consider the fatgraph with one positive number per every edge, producing the {\it metric fatgraph}. Then we attach the infinite stripe to the edge, with the width being the corresponding parameter. The transition functions rise from gluing stripes corresponding to edges into neighborhoods $U_v$, with the width being a positive parameter assigned to the edge. The key ideas of this description, which is due to Kontsevich \cite{konts} and further elaborated by Mulase and Penkava \cite{mulase}, lies within the theory of Strebel differentials. These are holomorphic quadratic differentials on a punctured surface with certain extra conditions. One can reconstruct the metric fatgraph and the corresponding complex structure for every Strebel differential so that their zeroes define the vertices of the fatgraphs, and the order of zero determines the valence of the corresponding vertex. At the same time, the punctures correspond to their double poles. All this can be summarized in the fact that Strebel differentials parametrize the trivial  
  $\mathbb{R}_+^s$-bundle over the moduli space of Riemann surfaces  with $s$ punctures.   

In Section 5, we use this fatgraph description to characterize the moduli space of $(1|1)$-supermanifolds with punctures: we use the term ``puncture" for marked points or $(0|1)$-divisors assigned to marked points on the underlying Riemann surface. At first, we consider the split $(1|1)$-supermanifolds, which can be viewed as Riemann surfaces with line bundle $\mathcal{L}$ over it. The corresponding moduli space can be then described by the flat $U(1)$ connections on the corresponding metric fatgraphs with zero monodromies around the punctures,  accompanied by a fixed divisor at punctures, one for every degree.   

Next, we describe this construction's deformation by expressing the tangent bundle's odd parts to the corresponding moduli space as \v{C}ech cocycles on the  Riemann surface $F$. These cocycles lead to the infinitesimal deformations of the transition functions, which could be continued beyond the infinitesimal level.

Parametrizing such \v{C}ech cocycles is a nontrivial problem, which, however, can be solved in the case when ${\rm deg}(\mathcal{L})=1-g-n-r/2$, where $n$ is the number of point punctures and $r$ is the (even) number of $(0|1)$-divisor punctures, and $g$ is a genus. In this case, the corresponding cocycles can be characterized by the ordered sets of complex odd parameters for every vertex, where the number of parameters in each set depends on the valence of the vertex. This is roughly twice more parameters than needed, so there are equivalences between complex structures constructed in such a way. We characterize those equivalences explicitly using sections of the appropriate line bundles.

Thus the fatgraph description of the split case, together with the para\-met\-rization of cocycles, immediately leads to complete parametrization of the complex structures of $(1|1)$-super\-manifolds with such degree. 

We note that on the level of uniformization, this is an important subclass of supermanifolds obtained in \cite{IPZ}, corresponding to flat connections with zero monodromies around punctures.

In Section 6, we use the results of Dolgikh, Rosly, and Schwarz \cite{schwarz}, who explicitly described the equivalence between $N=2$ super Riemann surfaces and $(1|1)$-supermanifolds, expressing the transition functions of $N=2$ SRS using the transition functions for $(1|1)$- supermanifolds obtained in Section 5.

In Section 7, we first discuss the involution on the moduli space of $N=2$ SRS, such that the fixed points of this involution are $N=1$ SRS. We then describe the split case, characterizing various 
choices of the corresponding line bundle using spin structures on the fatgraph, thus looking at the corresponding supermoduli space with the given assignment of R and NS punctures as a $2^{2g}$ covering space over moduli space of punctured Riemann surfaces. 
We then apply the involution to the deformations, first on the infinitesimal level and then continuing beyond, using the superconformal condition.
This eventually leads to our main Theorem \ref{maintheorem} which describes deformations of 
$N=1$ SRS.

\subsection*{Acknowledgements}
A.M.Z. is partially supported by Simons
Collaboration Grant 578501 and NSF grant
DMS-2203823.

\section{$(1|1)$-supermanifolds, $N=1$ and $N=2$ Super-Riemann Surfaces and Superconformal Transformations}

\subsection{Super Riemann surfaces and superconformal transformations.}
We remind that a {\it complex supermanifold} of dimension $(1|1)$  (see, e.g.,  \cite{br}) over some Grassmann algebra $S$ is a pair $(X,\mathcal{O}_{X})$, where $X$ is a topological space and $\mathcal{O}_X$ is a sheaf of supercommutative $S$-algebras over $X$ such that $(X,\mathcal{O}^{\rm{red}}_{X})$ can be identified with a Riemann surface (where $\mathcal{O}^{\rm{red}}_{X}$ is obtained from $\mathcal{O}_X$ by quoting out nilpotents) and for some open sets $U_{\alpha}\subset X$ and some linearly independent elements $\{\theta_{\alpha}\}$ we have $\mathcal{O}_{U_{\alpha}}=\mathcal{O}^{\rm{red}}_{U_{\alpha}}\otimes S[\theta_{\alpha}]$. 
We will also refer to $(X,\mathcal{O}^{\rm{red}}_{X})$ as a {\it base manifold}. These open sets $U_{\alpha}$ serve as coordinate neighborhoods for supermanifolds with coordinates $(z_{\alpha}, \theta_{\alpha})$.  The coordinate transformations on the overlaps $U_{\alpha}\cup U_{\beta}$ 
are given by the following formulas $z_{\alpha}=f_{\alpha\beta}(z_{\beta}, \theta_{\beta})$, $\theta_{\alpha}=\psi_{\alpha\beta}(z_{\beta}, \theta_{\beta})$, where $f_{\alpha\beta}$, $\psi_{\alpha\beta}$ are even and odd functions correspondingly. 
A super Riemann surface (SRS) $\Sigma$ \cite{CR},\cite{superwitten} over some Grassmann algebra $S$ is a complex supermanifold of dimension $1|1$ over $S$, with one more extra structure: there is 
an odd subbundle $\mathcal{D}$ of $T\Sigma$ of dimension $0|1$, such that for any nonzero section $D$ of $\mathcal{D}$ on an 
open subset $U$ of $\si$, $D^2$ is nowhere proportional to $D$, i.e. we have the exact sequence:
\begin{eqnarray}\label{exact}
0\to \mathcal{D}\to T\si\to \mathcal{D}^2\to 0.
\end{eqnarray}
One can pick the holomorphic local coordinates in such a way that this odd vector field 
will have the form $f(z,\theta)D_{\theta}$, where $f(z,\theta)$ is a non vanishing function and:
\begin{eqnarray}
D_{\theta}=\partial_{\theta}+\theta\partial_z, \quad D_{\theta}^2=\partial_z.
\end{eqnarray}
Such coordinates are called $superconformal$. The transformation between two superconformal coordinate systems 
$(z, \theta)$, $(z', \theta')$ is determined by the condition that $\mathcal{D}$ should be preserved, namely:
\begin{eqnarray}
D_{\theta}=(D_{\theta}\theta') D_{\theta'},
\end{eqnarray}
Locally one obtains:
\begin{eqnarray}
z'=u(z)+\theta\eta(z)\sqrt{\partial_zu(z)}, \quad \theta'=\eta(z)+\theta \sqrt{\partial_zu(z)+\eta(z)\partial_z\eta(z)},
\end{eqnarray}
so that the constraint on the transformation emerging from the local change of coordinates is  $D_{\theta} z'-\theta'D_\theta \theta'=0$.
%An important example of a super Riemann surface is the Riemann super sphere $SC^*$: there are two  
%charts $(z, \theta)$, $(z, \theta')$ so that
%\begin{eqnarray}
%z'=-\frac{1}{z},\quad \theta'=\frac{\theta}{z}.
%\end{eqnarray}

\subsection{$N=2$ super Riemann surfaces.} $N=2$ super Riemann surfaces ($N=2$ SRS)  is a generalization of super Riemann surfaces, being a supermanifold of dimension $(1|2)$ with extra structure.  Its tangent bundle has two subbundles $\mathcal{D}_+$ and $\mathcal{D}_-$, so that each of them are integrable, meaning that if $D_{\pm}$ are nonvanishing sections of $\mathcal{D}_{\pm}$, we have 
\begin{eqnarray}
D_{+}^2=aD_{+}, \quad D_{-}^2=bD_{-}
\end{eqnarray}
for some functions $a$ and $b$. At the same time, the direct sum $\mathcal{D}_+ \oplus \mathcal{D}_-$ is non-integrable, so that $[D_{+}, D_{-}]$ is a basis for the tangent bundle. Namely, for $N=2$ super Riemann surface $\Sigma$ one has the following exact sequence:
\begin{eqnarray}
0\to \mathcal{D}_+ \oplus \mathcal{D}_-\to T\Sigma\to \mathcal{D}_+ \otimes \mathcal{D}_-\to 0.
\end{eqnarray}
As in the case of super Riemann surfaces one can show that there exist superconformal coordinates in which locally $\mathcal{D}_{+}$ and $\mathcal{D}_{-}$ are generated by:
\begin{eqnarray}
D_{+}=\partial_{\theta_{+}}+\frac{1}{2}\theta_-\partial_z, \quad D_{-}=\partial_{\theta_{-}}+\frac{1}{2}\theta_+\partial_z,
\end{eqnarray}
so that $D_{\pm}^2=0$, $[D_{+}, D_{-}]=\partial_{z}$.

It turns out, there is an equivalence between $(1|1)$ supermanifolds and $N=2$ SRS as it was estabished by Dolgikh, Rosly and Schwarz \cite{schwarz}.  
One can notice that there is an involution $\theta_{+}\leftrightarrow \theta_{-}$. 
The corresponding complex $(1|1)$ supermanifold constructed from the $N=2$ SRS after the involution is of course generally a different one and it is called {\it dual}. In fact, such a dual supermanifold turns out to be a supermanifold of $(0|1)$ divisors of the original one. The self-dual $(1|1)$ supermanifolds are of course $N=1$ super Riemann surfaces.

We will discuss these questions in more detail later in the text.

\subsection{Punctures: Ramond and Neveu-Schwarz}
Let us now discuss the types of punctures on $N=1$ super Riemann surface.

The \emph{NS puncture} is a natural generalization of the puncture of ordinary Riemann surface, and can be considered as any point $(z_0,\theta_0)$ on the super Riemann surface. Locally one can associate to it a $(0|1)$-dimensional divisor of the form $z=z_0-\theta_0\theta$, which is the orbit with respect to the action of the group generated by $D$, and this divisor uniquely determine the point $(z_0,\theta_0)$ due to the superconformal structure.

Let us consider the case when the puncture is at $(0,0)$ locally. In its neighborhood let us pick a coordinate transformation 
\Eq{
z=e^{w}, \tab \theta= e^{w/2}\eta, \label{trans1}} such that the neighborhood (without the puncture) is mapped to a \emph{supertube} with $w$ sitting on a cylinder $w \sim w+2\pi i$, and $D_\theta$ becomes
\Eq{
D_{\theta}=e^{-w/2}(\partial_{\eta}+\eta\partial_{w}).
} Hence $(w,\eta)$ are superconformal coordinates, and we have the full equivalence relation given by \Eq{w \sim w+2\pi i, \tab \eta\to -\eta.\label{eq1}}

The case of \emph{Ramond puncture} is a whole different story. On the level of super Riemann surfaces, the associated divisor is determined as follows. In this case, we are looking at the 
case when the condition that  $D^2$ is linearly independent of $D$ is violated along some $(0|1)$ divisor. Namely, in some local coordinates $(z,\theta)$ near the Ramond puncture with coordinates $(0,0)$, $\mathcal{D}$ has a section of 
the form 
$$D^*_{\theta}=\partial_\theta+z\theta \partial_z.$$
We see that its square vanishes along the $Ramond$ $divisor$ $z=0$. One can map the neighborhood patch to the supertube using a different coordinate transformation 
\Eq{z=e^{w},\tab \theta=\eta\label{trans2},} those coordinates on the supertube will be superconformal, since \Eq{
D_{\eta}=\partial_{\eta}+\eta\partial_{w}.} Notice that the identifications we have to impose on $(w, \eta)$ now become: \Eq{w \sim w+2\pi i,\tab \eta\to+\eta.\label{eq2}}

To describe Ramond punctures globally, consider a subbundle $\mathcal{D}$ is generated by  such operators $D^*_{\eta}$, for the Ramond punctures $p_1, p_2,\dots p_{n_R}$. We have exact sequence:
\begin{eqnarray}
0\to \mathcal{D}\to T\Sigma\to \mathcal{D}^2\otimes \mathcal{O}(\mathscr{P})\to 0,
\end{eqnarray}
where $\mathscr{P}=\sum^{n_R}_{i=1}{\mathscr{P}_i}$ is a divisor where $\mathcal{D}^2=0$ mod $\mathcal{D}$. 

 In the split case $T\Sigma|_X=TX\oplus  {\mathcal{G}}$, dividing the tangent space to $T\Sigma$ into even and odd parts, which one can identify with $\mathcal{D}^2\otimes \mathcal{O}(\mathscr{P})$and $\mathcal{D}$ correspondingly. Also, notice that after reducing it to the base manifold $\mathcal{O}(\mathscr{P})=\mathcal{O}(\sum^{n_R}_{i=1}{p_i})$. 
Therefore,
$$\mathcal{G}^2=TX\otimes\mathcal{O}\Big(-\sum^{n_R}_{i=1}{p_i}\Big). $$
Since ${\rm deg}(TX)=2-2g$, that automatically implies that ${\rm deg}(\mathcal{G})=1-g-n_R/2$, leading to the fact that there should be even number of such punctures, known as Ramond or simply R punctures. We refer to the Section 4.2.2. of \cite{superwitten} for more details.

\section{Fatgraphs and spin structures}
From now on we will consider Riemann surfaces of genus $g$ with $s$ punctures  ($s>0$) and negative Euler characteristic, which we will denote as  $F^s_g$ or simply $F$. The corresponding closed version will be denoted as $F^c$.  

Consider the {\it fatgraph} $\tau$, corresponding to an $s$-punctured surface $F$.  
This is a graph, which is homotopically equivalent to $F$, with cyclic orderings on half-edges for every vertex \cite{pbook} induced by the orientation of the surface. 

 Let $\t_0,\t_1$ denote the set of vertices and edges of $\t$ respectively. Let $\w$ be an orientation on the edges $\t_1\subset\t$. As in \cite{penzeit}, we define a \emph{fatgraph reflection} at a vertex $v$ of $(\t,\w)$ to reverse the orientations of $\w$ on every edge of $\t$ incident to $v$.
\begin{Def}\label{fatgraph} We define $\cO(\t)$ to be the equivalence classes of orientations on a trivalent fatgraph $\tau$ spine of $F$, where the equivalence relation is given by $\w_1\sim \w_2$ iff $\omega_1$ and $\omega_2$ differ by finite number of fatgraph reflections. It is an affine $H^1$-space where the cohomology group $H^1:=H^1(F;\Z_2)$ acts on $\cO(\t)$ by changing the orientation of the edges along cycles.
\end{Def}

In \cite{penzeit, IPZ, huang2019super} various realizations of the spin structures on the surface $F$, characterized by a trivalent fatgraph $\tau$ are described. Those results can be easily generalized to a fatgraph with vertices of any valence.

In fact, following \cite{johnson1980spin}, a spin structure can be characterized by a quadratic form $q:H_1(F;\Z_2)\to \Z_2$ 
so that for any cycles $a,b$ one has $q(a+b)=q(a)+q(b)+a\cdot b$ where $a\cdot b$ denotes the intersection form. 

The space of all orientation classes is an affine $H^1$-space. Really, fix a fatgraph $\t$, we denote by $o_{\w}(e)$ the orientation of the edge $e\in\t$ in the orientation $\w$. We define $\d_{\w_1,\w_2}:\t_1\to \Z_2$ by 
\begin{equation}
\d_{\w_1,\w_2}(e):=\case{+1& o_{\w_1}(e)=o_{\w_2}(e),\\-1& o_{\w_1}(e)\neq o_{\w_2}(e),}
\end{equation}
which defines an element in $H^1(F;\Z_2)$.

\begin{Prop}\label{isospin}\cite{penzeit, IPZ} The set of spin structures is isomorphic to the space of quadratic forms $\cQ(F)$ on $H_1(F;\Z_2)$, and also isomorphic to $\cO(\t)$, as affine $H^1$-spaces. 
\end{Prop}

This leads to the following important consequence.

\begin{Thm}\cite{IPZ} Given an oriented simple cycle $\c\in\pi_1(F)$ homotopic to a path on the fatgraph $\t$ with orientation class $[\w]$, the corresponding quadratic form is given by
\Eq{\label{ourquad}
q([\gamma])=(-1)^{L_\gamma}(-1)^{N_\gamma}=(-1)^{R_\gamma}(-1)^{\over[N]_\gamma},}
where $L_\c$ (resp. $R_\c$) is the number of left (resp. right) turns of $\gamma$ on the fatgraph $\t$, and $N_\c$ (resp. $\over[N]_\gamma$) is the number of edges of $\t$ such that $\gamma$ and $\w$ have the same (resp. opposite) orientation.
\end{Thm}

If we talk about the paths corresponding to the boundary cycles on the fatgraph, $q([\gamma])=(-1)^k$, where $k$ is a number of edges with orientation opposite to the canonical orientation of $\gamma$. 
One can identify them with R and NS punctures in uniformization picture \cite{IPZ2}, so that $k$ is even for R punctures and odd for NS punctures.

 In fact there is another way of thinking about the spin structures, using graph connections \cite{fock1998poisson},\cite{pbook},\cite{bourque}.

\begin{Def}\label{graphconnect}\cite{pbook} Let $G$ be a group. A \emph{$G$-graph connection} on $\t$ is the assignment $g_e\in G$ to each oriented edge $e$ of $\t$ so that $g_{\over[e]}=g_e\inv$ if $\over[e]$ is the opposite orientation to $e$. Two assignments $\{g_e\},\{g_e'\}$ are equivalent iff there are $t_v\in G$ for each vertex $v$ of $\t$ such that $g_e'=t_v g_e t_w\inv$ for each oriented edge $e\in\t_1$ with initial point $v$ and terminal point $w$.
\end{Def}

Therefore, we obtain the following description of spin structures.

\begin{Cor}\cite{IPZ}
The space of spin structures on $F$ is identified with $\mathbb{Z}_2$-graph connections on a given fatgraph $\tau$ of $F$.
\end{Cor}
%\begin{proof}
%Pick an orientation $\omega$ on fatgraph $\tau$. For any orientation $\omega'$, we label the edge by $+1$ if they agree on that edge and by $-1$ if they do not agree. Classes of orientations give rise to classes of $\mathbb{Z}_2$-connections, since the reflection of orientation at the vertices corresponds to multiplication by $-1$.
%\end{proof}
We will refer to the fatgraph with the associated spin structure/$\mathbb{Z}_2$-graph connection as {\it spin fatgraph}.

\section{Complex structures and Strebel differentials}

\subsection{Gluing of Riemann surfaces.}
Consider the fatgraph $\tau$ corresponding to an $s$-punctured surface $F$ of genus $g$. 
Let us assign a positive real parameter $L_j$ associated to every edge $j$. We will refer to the resulting object as a {\it metric fatgraph}. 

It is known from the works of Penner (the so-called convex hull construction) that the  fatgraphs with valence of every vertex greater ore equal to 3 or dual ideal cell decompositions of Riemann surfaces describe the mapping class group-invariant cell decomposition of  the decorated Teichm\"uller space (see e.g. \cite{pbook}, \cite{penner}), a universal cover of $\mathbb{R}^s_+\otimes \mathcal{M}_{g,s}$, where $ \mathcal{M}_{g,s}$ is a moduli space of Riemann surfaces of genus $g$ with $s$ marked points. Then the trivalent fatgraphs correspond to the higher-dimensional cells of dimension $6g-6 +3s$. 

An important problem is how to reproduce inequivalent complex structures based on the data of metric fatgraphs.  In an important work of Mulase and Penkava \cite{mulase}, based on earlier ideas of Kontsevich \cite{konts}, allows to 
construct the appropriate covering of a Riemann surface and the transition functions, associated with a given fatgraph, thus exhausting all possible complex structures.
Let us have a look at those in detail. 

Fixing an orientation on $\tau$, we consider a neighborhood $U_v$ with coordinate $w$ corresponding to the fixed $m$-valent vertex $v$, so that the vertex is placed at the point $w=0$.
One can describe that neighborhood by considering stripes 
\begin{eqnarray}
\{z_j,\in \mathbb{C}, 0<{\rm Re} (z_i)<L_j\}, \quad j=1, \dots , m
\end{eqnarray}
glued together via formula
\begin{eqnarray}
w=e^{\frac{2\pi i(j-1)}{m}}z_j^{\frac{2}{m}}, \quad j=1, \dots , m
\end{eqnarray}
if all $m$ edges are pointing out from vertex $v$. In the case if one or more of them is pointing towards vertex $v$, we substitute the above formulas by 
\begin{eqnarray}
w=e^{\frac{2\pi i(j-1)}{m}}(L_j-z_j)^{\frac{2}{m}}, \quad j=1, \dots , m.
\end{eqnarray}
One can construct such coordinate patches around every such vertex. The overlaps $U_v\cap U_{v'}$ are described by the corresponding stripes associated to the edge $j$ of the fatgraph runnig between $v$ and $v'$. Note, that there is no triple intersections on such a punctured surface and that the vertices of the fatgraph belong to the boundary of the intersections.

Let us look at the transition functions on the overlaps between two such coordinate neighborhoods 
$U_v$, $U_{v'}$ around neighboring vertices  $v$ and $v'$, assuming the edge is pointing from $v$ to $v'$.  
We note that both coordinates $w$ and $w'$ are expressed in terms of $z_j$ in the following way:
\begin{eqnarray}
w=c_jz_{j}^{2/m}, \quad w'=c'_j(L_j-z_j)^{2/m},
\end{eqnarray}
where $c_j$, $c'_j$ are $m$th and $m'$th roots of unity. 
The resulting overlap coordinate transformation $f_{v'v}$  between patches is given by the following formula:
\begin{eqnarray}\label{overlap1}
w'=f_{v',v}(w)=c'_j(L_j-c_j^{-m/2} w^{m/2})^{2/m},
\end{eqnarray}
where ${-i\pi /2}<$ arg($w^{m/2}$)$ <{i\pi /2}$.  That completely describes the transition functions between charts for the punctured Riemann surface $F$.

Note that if the consecutive edges $L_1, L_2, \dots, L_n$ correspond to boundary piece of the fatgraph associated with puncture $p$, one can glue the following coordinate  neighborhood $V_p$ with coordinate $y$ covering the puncture:
\begin{eqnarray}
y= e^x=\exp\Big({\frac{2\pi i}{a_B}(L_1+\dots L_{k-1}+z_k)}\Big),~ {\rm where} ~a_B=L_1+\dots +L_n,
\end{eqnarray}
so that one glues together top or bottom part of the stripes based on orientation and $x$-variable is on a cylinder $x\sim x+2\pi i$. Suppose the $k$-th strip is glued to the vertex $v_k$ with coordinate $w_k$ as above, then
the transition function $f_{pv}$ is given by:
\begin{eqnarray}\label{overlap2}
y=f_{p,v}(w)=\exp\Big({\frac{2\pi i}{a_B}(L_1+\dots L_{k-1}+w_k^{m/2})}\Big).
\end{eqnarray}

 In the following we will adopt the following notation for $L$-parameters: if the edge connects two vertices $v,v'$ we will denote the corresponding parameter as $L_{v,v'}$.

\subsection{Strebel differentials.} An important object in the constructions of \cite{konts},\cite{mulase} are the {\it Strebel differentials}, the quadratic meromorphic differentials with special  properties. A nonzero quadratic differential is a holomorphic section $\mu$ of $K^{\otimes 2}$, where $K$ stands for canonical bundle on $F$. It  defines a flat metric on the complement of the discrete set of its zeroes, written  in local coordinates as $|\mu(z)|dzd\bar z$,  where $\mu=\mu(z)dz^2$. 

A {\it horizontal trajectory} of a quadratic differential is a curve along which $\mu(z)dz^2$ is real and positive. The {\it Strebel} differential is the one for which the union of nonclosed  trajectories has measure zero.
Non-closed trajectories of a given Strebel differential decompose the surface into the
maximal ring domains swept out by closed trajectories. These ring domains
can be annuli or punctured disks. All trajectories from any fixed maximal
ring domain have the same length, the circumference of domain. 

The following Theorem is due to Strebel:

\begin{Thm}\cite{strebel}\label{streb}
For any connected closed Riemann surface $F^c$ with s distinct points $p_1,\dots, p_s$, $s>0$ and genus $g$, $s>\chi(F^c)=2-2g$
and n positive real numbers $a_1, \dots, a_s$ there
exists a unique Strebel differential on $F=F^c\backslash\{ p_1,\dots, p_s\}$, whose maximal ring
domains are $s$ punctured disks surrounding $p_i$'s with circumference $a_i$'s.
\end{Thm}

The union of non-closed trajectories of Strebel differentials together with their zeroes define a graph, embedded into a Riemann surface, thus giving to it a fatgraph structure. Every vertex of a fatgraph corresponds to the zero of the  Strebel differential of degree $m-2$, where $m\ge 3$ is the valence
of the vertex. The length of each edge gives the graph a metric structure. 

For every such Strebel differential one can construct the covering associated to the corresponding fatgraph, described in the previous section and vice versa,  so that in the charts $U_v$, Strebel differential $\mu$ has the following explicit form:
\begin{eqnarray}
\mu |_{U_v}=w^{m-2}dw^2.
\end{eqnarray}
 It also has pole of order $2$ at punctures, so that in $y$-coordinates for each neighborhood $V_p$ the differential looks as follows: 
 \begin{eqnarray}
 \mu |_{V_{p}}=-\frac{a^2_p}{4\pi^2}y^{-2}dy^2.
 \end{eqnarray} 

 One can formulate then the following Theorem.
 
 \begin{Thm} \cite{konts}
 Let $\mathcal{M}^{\rm comb}_{g,s}$ denote the set of equivalence classes of connected
ribbon graphs with metric and with valency of each vertex greater than or
equal to 3, such that the corresponding noncompact surface has genus g and s
punctures numbered by $1,\dots, s$. The map $\mathcal{M}_{g,s}\times \mathbb{R}^s\to \mathcal{M}^{\rm comb}_{g,s}$
which associated to the surface $F^c$ and numbers $a_1, \dots a_s$ 
 the critical graph of the canonical Strebel differential from Theorem \ref{streb} is one-to-one. 
 \end{Thm}
 
 In this paper we do not need more properties of Strebel differentials, however, we refer reader to \cite{mulase}, as well as original source \cite{strebel} for more information.

\section{Complex structures on $(1|1)$ supermanifolds}

\subsection{Split case}

Let us consider the punctured Riemann surface glued as in previous subsection 
using metric fatgraph and overlapping neighborhoods ${U}_v$ corresponding to vertices. 
To construct the cooordinate transformations for a split $(1|1)$-supermanifold $SF$ with such a base complex manifold, one has to consider a line bundle $\mathcal{L}$ over $F^c$. Then the coordinate transformations for the coordinates  $(w', \xi')$, $(w, \xi)$, $(y,\eta)$ corresponding to  
neighborhoods $U_v, U_{v'}, V_p$ of vertices $v, v'$ and puncture $p$ are given by the following formulas:
\begin{eqnarray}\label{gsplit}
&&\xi'=g_{v',v}(w)\xi, \quad w'=f_{v',v}(w),\nonumber\\
&&\eta=g_{p,v}(w)\xi, \quad y=f_{p,v}(w),
\end{eqnarray} 
where $g_{v',v}$, $g_{p,w}$ is the holomorphic function, serving as a transition function
of bundle $\mathcal{L}$. 
The collection $\{g_{v',v}, g_{v,p}\}$ generates a \v{C}ech cocycles 
\begin{eqnarray}
g_{v,v'}|_{U_{v}\cap U_{v'}}\in H^0(U_{v}\cap U_{v'}, \mathcal{O}^*),\quad 
g_{v,p}|_{U_{v}\cap V_p}\in H^0(U_{v}\cap V_p, \mathcal{O}^*),
\end{eqnarray} 
representing the Picard group of $F^c$, i.e. $\check{H}^1(F, \mathcal{O}^*)$, if the following constraint on $g_{v',v}$ and $\{g_{v,p}\}$ is imposed around the given puncture $p$:
\begin{eqnarray}\label{bc}
g_{v,v'}|_{U_{v'}\cap U_{v}\cap V_p}=g_{v,p}g_{p,v'}.
\end{eqnarray} 

Then the following Proposition holds.

\begin{Prop}\label{deg0}{When $\mathcal{L}$ is degree 0 over $F^c$, the fatgraph data 
describing it is a $U(1)$ graph connection with a trivial monodromy around every boundary piece.
%ii) Given a bundle $\mathcal{F}$ of a fixed degree $d$ on $F^c$ with transition functions 
%$g^{(\alpha, \beta)}_{v',v}$, one can construct all transition functions for all other elements of the Picard group of degree $d$ by considering $\tilde{g}_{v',v}=e^{2\pi i h_{v'v}}g^{(\alpha, \beta)}_{v',v}$, where 
%elements $e^{2\pi i h_{v'v}}$ generate $U(1)$ graph connection with a trivial monodromy around every puncture
} 
\end{Prop}
\noindent {\bf Proof.}
 Notice, that one can choose $g_{v'v}$ to be constant functions with values on a unit circle, which on the level of fatgraph is described by $U(1)$-graph connection, so that  $g_{v,v'}=e^{ih_{v,v'}}$, where $h_{v,v'}\in \mathbb{R}$ .  Indeed the corrresponding holomorphic equivalences for the corresponding \v{C}ech cocyle reduce to constant $U(1)$ gauge transformations at the vertices. 
 However, according to the condition (\ref{bc}) that we imposed, we have to have $g_{v_1,v_2}g_{v_2,v_3}\dots g_{v_{n-1},v_{n}}g_{v_{n},v_{1}}=1$, which is exactly the trivial monodromy condition.
%To prove ii), notice that zero monodromy U(1)-connections enumerate all degree zero bundles on 
%$F^c$ $$ 
\hfill $\blacksquare$\\

In order to describe any line bundle of degree $d$, one has to do the following. 
First, choose a fixed divisor of degree $d$, say a linear combination of puncture points. Then, multiplying it on appropriate bundle of degree 0, one can reproduce the original bundle. 
Since we described the moduli spaces of degree 0 bundles in a Proposition \ref{deg0} above, we can now characterize split punctured supermanifolds.

\begin{Thm}
Consider the following data on a fatgraph $\tau$:

\begin{itemize}
\item Metric structure.\\

\item  Flat $U(1)$-connection with zero monodromies around punctures.\\

\item Fixed divisor $M$ of degree $d$, which is a linear combination of puncture points. 
\end{itemize}
The data above determines the complex split $(1|1)$-supermanifold corresponding to the line bundle of degree $d$ on $F$.  For a fixed divisor $M$, metric fatgraphs with $U(1)$ connections describe the moduli space of split $(1|1)$-supermanifolds.
\end{Thm}

\subsection{Infinitesimal Deformations and various types of punctures} 
As usual, the infinitesimal deformations of the above formulas leading to generic non-split structure are described by $H^1(SF^c, ST)$, where $ST$ is a tangent bundle of $SF^c$, where $SF^c$ is a split $(1|1)$-supermanifold which we discussed in the previous section. 
Since we are deforming the split case one can describe infinitesimal deformations $\rho\in H^1(SF^c, ST)$ using \v{C}ech cocycles, 
i.e. in coordinates $(\xi,w)$ of $U_v$ on $U_v\cap U_{v'}$ and $U_v\cap V_{p}$:
\begin{eqnarray}\label{superco}
&&\rho_{v ,v'}={\rm v}_{v,v'}(w)\partial_w+\xi\alpha_{v,v'}(w)\partial_w+
\beta_{v,v'}(w)\partial_{\xi}+{\rm u}_{v,v'}(w)\xi\partial_{\xi},\\
&&\rho_{p ,v}={\rm v}_{p,v}(w)\partial_w+\xi\alpha_{p,v}(w)\partial_w+
\beta_{p,v}(w)\partial_{\xi}+{\rm u}_{p,v}(w)\xi\partial_{\xi}.\nonumber
\end{eqnarray}
where the indices $v,v'$ and $p,v$ here mean that the corresponding elements are the corresponding \v{C}ech cocycles considered on the intersections $U_v\cap U_{v'}$, and $V_p\cap U_{v}$. 

Now we need to specify the behavior at the punctures to describe the cocycles $\rho$ leading to deformations of $SF$ in terms of cocyles on $F^c$.

There are two types of punctures  we want to consider:\\
\begin{itemize}
\item puncture as a $(0|1)$-dimensional divisor on $SF^c$. We denote the number of such punctures as $r$.\\
\item puncture as a $(0|0)$-dimensional divisor, or in other words just a point on $SF^c$. 
We denote the number of such punctures as $n$.\\
 \end{itemize}

 Let $T$ be the tangent bundle of $F^c$, $D_{n+r}$ be the divisor corresponding to sum of all points on $F^c$ corresponding to punctures, and $D_n$ is the sum of the ones corresponding to point-punctures on $SF$. Let us look in detail at the components of (\ref{superco}): 
\begin{eqnarray}
&&{\rm v}\in \check{Z}^1( F^c, T\otimes \mathcal{O}(-D_{n+r})), \quad {\rm u}\in \check{Z}^1(F^c, \mathcal{O}), \\
&&\beta\in \Pi\check{Z}^1(F^c, {\mathcal{L}}\otimes \mathcal{O}(-D_{n})),\quad 
\alpha\in \Pi\check{Z}^1(F^c,T\otimes {\mathcal{L}}^{-1}\otimes \mathcal{O}(-D_{n+r})),\nonumber
\end{eqnarray}
where $\check{Z}^1$ is the notation for \v{C}ech cocycles of degree 1. Note, that 
we need to impose the constraints on cocycles on $V_p\cap U_v\cap U_{v'}$: 
\begin{eqnarray}
s_{v,v'}|_{U_{v'}\cap U_{v}\cap V_p}=s_{v,p}+s_{p,v'},
\end{eqnarray} 
where $s={\rm v,u}, \alpha, \beta$. 
Here  u- and v-terms correspond to the deformations of the original manifold $F$, and 
notice, that we already incorporated the moduli for the base manifold $F$ and the line bundle $\mathcal{L}$ in the formulas (\ref{gsplit}). The odd deformations, provided by the cycles  $\alpha,\beta$ give the following deformations for the upper line of (\ref{gsplit})
\begin{eqnarray}
\xi'=g_{v'v}(w)(\xi+\beta_{v',v}(w)), \quad w'=f_{v'v}(w+\xi \alpha_{v',v}(w)),
\end{eqnarray}
which describes (in the first order in complex parameters)all possible complex structures on the punctured supermanifold. 

If we remove the infinitesimality condition, formulas above will be deformed. Let us formulate it in a precise form.

\begin{Thm}\label{sman}
\begin{enumerate}
\item Consider the following data:\\

\begin{itemize}
\item A metric fatgraph with a $U(1)$-connection with trivial monodromy around boundary pieces, a fixed divisor which is a iinear combination of puncture points of degree $d$, which defines a split punctured $(1|1)$ supermanifold determined by base Riemann surface $F$ and line bundle $\mathcal{L}$.\\

\item  \v{C}ech cocycles $$\tilde{\beta}=\sum_k\sigma_k^{\beta}b_k, ~\tilde{\alpha}=\sum_k
\sigma_k^{\alpha}a_k,$$ so that 
$\{\sigma_k^{\beta}\}$, $\{\sigma_k^{\alpha}\}$ are two sets of odd parameters, 
\begin{eqnarray}
&&b_k\in \check{Z}^1(F^c, {\mathcal{L}}\otimes \mathcal{O}(-D_{n})), \nonumber\\ 
&&a_k\in \check{Z}^1(F^c,T\otimes {\mathcal{L}}^{-1}\otimes \mathcal{O}(-D_{n+r})),\nonumber 
\end{eqnarray}
where $D_{n+r}$ is the divisor corresponding to sum of all $s=n+r$ punctures on the closed surface $F^c$, $D_n$ is the sum of the certain subset of the set of punctures, and the cohomology classes of $\{b_k\}$  $\{a_k\}$ form a basis in the corresponding cohomology spaces.
\\
\end{itemize}

This data gives rise to a family of complex structures on $SF$, the $(1|1)$-supermanifold with $n$ point punctures and r $(0|1)$-divisor punctures, so that the transition functions on 
$SF$  are given by the following formulas on the overlaps $\{U_v\cap U_{v'}\}$:
\begin{eqnarray}\label{gnsplit}
\xi'=g^{(\alpha, \beta)}_{v'v}(w)(\xi+\beta_{v',v}(w)),\quad  w'=f^{(\alpha, \beta)}_{v'v}(w+\xi \alpha_{v',v}(w)),
\end{eqnarray}
where $g^{(\alpha, \beta)}_{v',v}, f^{(\alpha, \beta)}_{v',v},$ are holomorphic functions on the overlaps, depending on the parameters $\sigma_k^\alpha$ and $\sigma_k^\beta$ such that:
$$g^{(0,0)}_{v',v}=g_{v',v}, ~ f^{(0,0)}_{v',v}=f_{v',v},$$ where $\{f_{v'v}\}$, $\{g_{v',v}\}$ define the split supermanifold with the line bundle $\mathcal{L}$ and $s$ punctures so that in the first order in  $\{\sigma_k^\alpha\}$ and $\{\sigma_k^\beta\}$ we have 
$$
\tilde{\beta}_{v',v}(w)=\beta_{v',v}(w), ~\tilde{\alpha}_{v',v}(w)=\alpha_{v',v}(w).
$$

\item Let us fix the choice of transition functions in (\ref{gnsplit}), for every metric fatgraph $\tau$ with the U(1)-connection, divisor of degree d, and the odd data given by the cocycles $\tilde \beta, \tilde\alpha$ on $F^c$. 

The complex structures constructed in such a way are inequivalent to each other, and the set of such complex structures constructed by varying $\tau$ and the data on it, form a dense subset of maximal dimension in the moduli space of punctured $(1|1)$ supermanifolds with underlying line bundles of degree $d$. 
%Moreover, the cycles cohomologically equivalent to $\tilde{\beta}$, $\tilde{\alpha}$ produce a complex structure, equivalent to the one corresponding to 
%$\tilde{\beta}$, $\tilde{\alpha}$. 
\end{enumerate}
\end{Thm}
\noindent {\bf Proof.} 
 Let us look at the formulas (\ref{gnsplit}) as generic ones, for arbitrary holomorphic  functions $\{\alpha_{v'v}\}$, $\{\beta_{v'v}\}$ on overlaps. There is a finite number of odd parameters which parametrize all $\{{\alpha_{v',v}}\}$, $\{\beta_{v',v}\}$ corresponding to inequivalent complex structures. Expanding the formulas (\ref{gnsplit}) in terms of these parameters we obtain that in the linear 
order $\beta\in \Pi\check{Z}^1(F^c, {\mathcal{L}}\otimes \mathcal{O}(-D_{n}))$ and 
$\alpha\in \Pi\check{Z}^1(F^c,T\otimes {\mathcal{L}}^{-1}\otimes \mathcal{O}(-D_{n+r}))$ as in the infintesimal case. 
Conversely, since $\alpha$, $\beta$ represent the tangent space to the moduli space 
of complex structures, parameters $\sigma^{\alpha}$, $\sigma^{\beta}$ serve as 
coordinates there. Considering the corresponding 1-parametric subgroups generated by $\tilde \alpha$, $\tilde \beta$ we obtain formulas from (\ref{gnsplit}). 
The fact that the cohomologically equivalent cocycles lead to the equivalent complex structures is justified by dimensional reasons.
\hfill $\blacksquare $\\

It is, however, nontrivial to explicitly parametrize those cocycles $\alpha, \beta$. In the next subsection we will analyze the special case of supermanifolds with the line bundle $\mathcal{L}$ of negative degree.

\subsection{$(1|1)$-supermanifolds  with  ${\rm deg }(\mathcal{L})=1-g-r/2$}
 It is not easy to explicitly parametrize cocycles $\alpha, \beta$ from a fatgraph data if one does not fix a degree. From now on, we will be interested in the case when  ${\rm deg}(\mathcal{L})=1-g-k$, where $s\ge k\ge 0$ on $F^c$. Assuming that the number of divisor punctures is even and setting $k=r/2$, both bundles $\mathcal{L}\otimes \mathcal{O}(-D_{n})$ and $\mathcal{L}^{-1}\otimes T\otimes \mathcal{O}(-D_{n+r})$ have equal degree $1-g-r/2-n$ on $F^c$.

Let us be generic enough first and characterize the cycles in $\Pi\check{Z}^1(F^c, \mathcal{L}\otimes \mathcal{O}(-D_{n}))$ where $s\ge k=g-1-\deg{\mathcal{L}}\ge 0$,  using the data from the fatgraph. To do that, we define a cocycle $\rho$, a representative of $\Pi\check{H}^1(F^c, \mathcal{L}\otimes \mathcal{O}(-D_{n}))$ as follows: 
\begin{eqnarray}\label{cycles}
&&\rho_{v,v'}|_{U_v\cap U_{v'}}=\rho_v-\rho_{v'},~{\rm  so~ that ~ } ~\rho_v|_{U_v}=\frac{\sigma_v(w)}{w^{m_v-2}},~  \rho_{v'}|_{U_v'}=\frac{\sigma_{v'}(w')}{w'^{~m_{v'}-2}},\nonumber \\
&&\rho_{v,p}|_{U_v\cap V_p}=\rho_v ,
%+\lambda_p,~ {\rm so~ that }~ \lambda_p|_{V_p}=y\nu_p  
%-\rho_p,  ~{\rm  so~ that ~ } \rho_{p}|_{V_p}\in \check{H}^0(V_p, K^{1/2}).
%y^2\sigma_p dy^{1/2}
%{\rm ~and} ~ \rho_v~ {\rm is ~defined ~on~} V_p.
\end{eqnarray}
where $\rho_v$, $\rho_{v'}$ are meromorphic sections of $\mathcal{L}\otimes \mathcal{O}(-D_{n})$ on $U_v$, $U_{v'}$ correspondingly,  so that $m_v$ is valence of the given vertex $v$,   
$$\sigma_v(w)=\sum^{m_v-3}_{i=0}\sigma^i_v w^i$$ are the polynomials with odd coefficients, assigned to each fatgraph vertex $v$ of degree at most $m_v-3$,
Let us denote for simplicity $\tilde{\mathcal{L}}=\mathcal{L}\otimes \mathcal{O}(-D_{n})$.

%and $\nu_p$ is an odd number associated to every puncture $p$

Then the following proposition holds. 

\begin{Thm} \label{cycdes}

\begin{enumerate}
\item The cycles (\ref{cycles}) are uniquely defined by the numbers $\sigma_v$ at the fatgraph vertices, thus forming a complex vector space of dimension $4g-4+2s$. \\
\item  Cycle $\rho$ is cohomologous to cycle $\tilde \rho$ in $\Pi\check{H}^1(F^c, \tilde{\mathcal{L}})$ if and only if 
\begin{eqnarray}
\sigma_v(w)-\tilde{\sigma}_v(w)=\gamma^{(m-3)}(w), 
\end{eqnarray}
for every vertex $v$,  where $\gamma\in \Pi{H}^0(F^c, \tilde{\mathcal{L}}\otimes K^2\otimes \mathcal{O}(2D_{n+r}))$, so that $\gamma |_{U_v}=\gamma(w)$, $\gamma^{(m-3)}(w)$ is the Taylor expansion of $\gamma(w)$ up to order $m-3$.\\
\item  The cohomology classes of cycles $\rho$ span $\Pi\check{H}^1(F^c, \tilde{\mathcal{L}})$.
\end{enumerate}
\end{Thm}
\noindent {\bf Proof.} To prove part $(1)$ one just has to count the number of vertices and parametres at verrtices. An elememtary Euler characteristic computation shows that 
\begin{eqnarray}
2g-2+s=\sum_{j\ge 3}\Big(\frac{j}{2}-1\Big)\mathcal{V}_j(\tau),
\end{eqnarray}
where $\mathcal{V}_j(\tau)$ is the number of  $j$-valent vertices in $\tau$. Notice, that for a $j$-valent vertex $v$, we have exactly $j-2$ odd parameters from the expansion of $\sigma_v(w)$, which immediately leads to the necessary parameter count, giving $4g-4+2s$.

To prove $(2)$,  on each coordinate neighborhood $U_v$, Strebel differential 
$\mu$  has the form $\mu|_{U_v}=w^{m_v-2}dw^2$, and $\mu |_{V_p}=-\frac{a^2_p}{4\pi^2}\frac{dy^2}{y^2}$, which means that one can rewrite the formula for the cocycle 
\begin{eqnarray}
\rho_{v,v'}=(\gamma_v-\gamma_{v'})/{\mu}, \quad \rho_{v,p}=(\gamma_v-\gamma_{p})/{\mu},
\end{eqnarray}
where $\gamma_v|_{U_v}={\sigma_v}(w)$, $\gamma_p|_{V_p}=0$, so that $\gamma_v\in \Pi\check{H}^0(U_v, \tilde{\mathcal{L}}\otimes K^2)$ and $\gamma_p=0\in \Pi\check{H}^0(V_p, \tilde{\mathcal{L}}\otimes K^2)$ .

Suppose that such cocycle is exact, namely: 
\begin{eqnarray}
 \gamma_v/\mu-\gamma_{v'}/\mu=(a_v-a_{v'})|_{U_v\cap U_v'}, \quad  \gamma_v/\mu-\gamma_{p}/\mu=(a_v-a_{p})|_{U_v\cap U_p}
\end{eqnarray}
for all $v$ and $v'$, so that $a_v\in \Pi\check{H}^0(U_v, \tilde{\mathcal{L}})$, $a_p\in \Pi\check{H}^0(V_p, \tilde{\mathcal{L}})$. 
 It is equivalent to $(\gamma_v-a_v\mu)=(\gamma_{v'}-a_{v'}\mu)|_{U_v\cap U_v'}$, 
 $(\gamma_v-a_v\mu)=(\gamma_{p}-a_{p}\mu)|_{U_v\cap V_p}$, i.e. formulas
$$\gamma_v-a_v\mu=\gamma |_{U_v}, \quad \gamma_p-a_p\mu=\gamma |_{V_p}$$ 
defines  $\gamma$ as a holomorphic section on $F$, i.e. $\gamma\in \Pi\check{H}^0(F, \tilde{\mathcal{L}}\otimes K^2)$. Assuming $\gamma |_{U_v}=\gamma(w)$ and $a_v|_{U_v}=a(w)$, the identity $\gamma_v=a_v\mu+\gamma |_{U_v}$ is only possible if 
$$a_v(w)=\frac{\gamma(w)-\gamma^{(m-3)}(w)}{w^{m-2}}~{\rm  and}~ \sigma_v(w)=\gamma^{(m-3)}(w), $$
where $\gamma^{(m-3)}(w)$ is the Taylor expansion of $\gamma(w)$ up to order $m-3$.
Also the identity $\gamma_p=a_p\mu+\gamma |_{V_p},$ i.e. 
$$a_p\mu+\gamma |_{V_p}=0$$ is possible only 
if $\gamma |_{V_p}$ has poles not greater than 2 at $y=0$, or, more precisely,  $\gamma\in H^0(V_p,\tilde{\mathcal{L}}\otimes K^2\otimes \mathcal{O}(2D_{n+r}))$. 
Therefore, cycles 
$\rho$ and $\tilde \rho$ are cohomologous to each other iff  the relation between the parameters on the fatgraph $\{\sigma_v\}$ and $\{\tilde{\sigma}_v\}$ correspodingly parametrizing them is as follows: 
\begin{equation}\label{sigma}
\tilde{\sigma}_v(w)-\sigma_v(w)=\gamma^{(m-3)}(w)
\end{equation}
where  $\gamma\in H^0(F^c,\tilde{\mathcal{L}}\otimes K^2\otimes  \mathcal{O}(2D_{n+r}))$ on $F^c$, such that 
$\gamma |_{U_v}=\gamma(w)$ with poles at the punctures of $F$ of order less or equal to 2 so that $\gamma^{(m-3)}(w)$ is the Taylor expansion of $\beta(w)$ up to order $m-3$.

Now, to prove part $(3)$, we need to show that such classes of cocycles form a $2g-2+k+n$-dimensional complex space as elements of $\Pi\check{H}^1(F^c, \tilde{\mathcal{L}})$. For a given section $\gamma$ of $\tilde{\mathcal{L}}\otimes K^2\otimes \mathcal{O}(2D_{n+r})$, the collection of the coefficients in $\gamma^{(m-3)}$, for each  vertex $v$ form a vector in our $4g-4+2s$-dimensional space of $\sigma$-parameters. The space, spanned by all such vectors is a complex $2g-2+2s-k-n$-dimensional space.  Indeed, it cannot be of smaller dimension, since we know that $dim_{\mathbb{C}}\check{H}^1(F^c, \mathcal{L}\otimes \mathcal{O}(-D_{n}))=2g-2+k+n$, at the same time it cannot be of greater dimension, since we know that the dimension of space of such meromorphic global sections of $\tilde{\mathcal{L}}\otimes K^2\otimes \mathcal{O}(2D_{n+r})$  is $2g-2+2s-k-n$ by the Riemann-Roch theorem.  \hfill $\blacksquare$\\

%If $\lambda_p$ are non-trivial, cocycles  

%It is easy to see that cocycles $\rho$ are not homologous to each other on a punctured surface. It means, that in order to get a 

Now we are ready to formulate a Theorem regarding parametrization of complex structures via fatgraph data. 

\begin{Thm}
\begin{enumerate} 
\item Consider the following data associated to the fatgraph $\tau$:\\

\begin{itemize}
\item Metric structure and  a $U(1)$-connection on $\tau$ with zero monodromy around punctures and a fixed divisor of degree $d=1-g-r/2$ at the punctures.\\

\item Two complex odd parameter sets $\{\sigma^{\alpha}_{v,k}\}$, $\{\sigma^{\beta}_{v,k}\}$ at each vertex $v$, so that $k=0,\dots, m_v-3$. 
\\
\end{itemize}
We will call two sets of  data from {\rm (1)} associated to fatgraph $\tau$ equivalent if the odd data are related as in Theorem \ref{cycdes}.\\

Constructing transition functions $f_{v',v}$ and $g_{v',v}$ from the even fatgraph data and cocycles 
$\tilde{\alpha}$, $\tilde{\beta}$, corresponding to $r$ $(0|1)$-divisor punctures and $n$ point punctures from the odd data, one obtains a family of complex structures on $(1|1)$-supermanifold in the framework of Theorem \ref{sman}. \\

\item Fixing the transition functions in (\ref{gnsplit}) and considering one such complex structure per equivalence class  of data for every fatgraph $\tau$, we obtain a set of inequivalent complex structures, which is a dense subspace of odd complex dimension $4g-4+2n+r$  in the space of all complex structures on $(1|1)$-supermanifolds with base line bundle of degree 
$d=1-g-r/2$ and $s=n+r$ punctures, where $n$ is the number of point punctures and r is the number of $(0|1)$-divisor punctures. 
\end{enumerate}
\end{Thm}
\noindent {\bf Proof.} The first part of data allows to construct split $(1|1)$-supermanifold as we know from previous sections, the odd data from the second allows to construct the corresponding cycles  $\tilde{\alpha}\in \Pi\check{Z}^1(F^c,\mathcal{L}\otimes \mathcal{O}(-D_{n}))$ and $\tilde{\beta}\in \Pi\check{Z}^1(F^c, \mathcal{L}^{-1}\otimes T\otimes \mathcal{O}(-D_{n+r}))$. If we choose an  orientation on the fatgraph, the formulas (see Theorem \ref{sman}):
\begin{eqnarray}\label{gnsplit2}
\xi'=g^{(\alpha, \beta)}_{v',v}(w)(\xi+\beta_{v',v}(w)),\quad  w'=f^{(\alpha, \beta)}_{v',v}(w+\xi \alpha_{v',v}(w))
\end{eqnarray}
produce the transition functions on $U_v\cap U_{v'}$ for the vertex oriented from $v$ to $v'$.
\hfill $\blacksquare$
\\

\noindent {\bf Remark.} Note, that the gauge equivalence for $U(1)$ connection, produce the following identification. If real numbers $h_{v,v'}$ parametrize $U(1)$ connection, then the transformations:
\begin{eqnarray}
&&h_{v,v'}\to h_{v,v'}+t_v-t_v', \nonumber \\
&&\sigma^{\alpha}_v, \sigma^{\beta}_v\to e^{it_v}\sigma^{\alpha}_v,  e^{-it_v}\sigma^{\beta}_v
\end{eqnarray}
produce equivalent configuration for infinitesimal parameters $\sigma$.
 In the paper \cite{IPZ} the uniformization version of $N=2$ Teichm\"uller space was constructed (see also \cite{baranov1989geometry}, \cite{natanzon2004moduli}), which corresponds exactly to $(1|1)$-supermanifolds, which serves as a universal cover for the 
 one we use here in the case of ${\rm deg}(\mathcal{L})=1-g-r/2$.   The above identifications played an instrumental role in the construction.\hfill $\blacksquare$\\

In the next two sections we will use the obtained results to describe transition functions for the corrresponding dual supermanifold and $N=2$ super Riemann surface following \cite{schwarz}.

\subsection{Dual $(1|1)$ supermanifold.} Finally, we give a desription of the concept dual $(1|1)$ supermanifold is a supermanifold of $(0|1)$ divisors of $SF$. To describe the explicit coordinates and coordinate transformations  transformations on such an object one can use a very simple equation (see, e.g., \cite{superwitten}):
\begin{equation}
w=a+\zeta\xi
\end{equation}
where $a,\zeta$ are the coordinates parametrizing such a $(0|1)$ divisor. Let us derive the formulas for the transformations of $a,\zeta$ variables, for the transformation between the charts with coordinates $(a, \zeta)$ and $(a',\zeta')$, so that $w'=a+\zeta'\xi'$.
\begin{eqnarray}
&&\xi'=g^{(\alpha, \beta)}_{v'v}(a+\zeta\xi)(\xi+\beta_{v',v}(a+\zeta\xi)),\nonumber\\ 
&& a'+\zeta'\xi'=f^{(\alpha, \beta)}_{v',v}(a+\zeta\xi+\xi \alpha_{v',v}(a+\zeta\xi)).\nonumber
\end{eqnarray}
We will now substitute first equation in the second and obtain:
\begin{eqnarray}
&& a'+\zeta'g^{(\alpha, \beta)}_{v'v}(a+\zeta\xi)(\xi+\beta_{v',v}(a+\zeta\xi))=\nonumber\\
&&f_{v'v}(a+\zeta\xi+\xi \alpha_{v',v}(a)).\nonumber
\end{eqnarray}
which leads to two equations:
\begin{eqnarray}\label{aprime}
&&a'+\zeta' g^{(\alpha, \beta)}_{v'v}(a)\beta_{v',v}(a)=f^{(\alpha, \beta)}_{v'v}(a)\\
&&\zeta'g^{(\alpha, \beta)}_{v'v}(a)+\zeta'\zeta\partial_a (g^{(\alpha, \beta)}_{v'v}(a)\beta_{v',v}(a))=\zeta\partial_a f^{(\alpha, \beta)}_{v'v}(a)-\partial_a f^{(\alpha, \beta)}_{v'v}(a)\alpha_{v',v}(a)\nonumber
\end{eqnarray}
The latter equation  immediately gives the transformation for $\zeta$:
\begin{eqnarray}
\zeta'=g^{(\alpha, \beta)}_{vv'}(a)(1+\zeta g^{(\alpha, \beta)}_{vv'}(a)\partial_a (g^{(\alpha, \beta)}_{v'v}(a)\beta_{v',v}(a)))(\partial_af_{v'v}(a)\zeta- \partial_af_{v'v}(a)\alpha_{v',v}(a))\nonumber
\end{eqnarray}
which could be simplified as follows:
\begin{eqnarray}
&&\zeta'=g^{(\alpha, \beta)}_{v,v'}(a)\partial_a f^{(\alpha, \beta)}_{v'v}(a)\Big(s_{v,v'}(a)\zeta- \alpha_{v',v}(a)\Big)\\
&&s_{v,v'}=\Big(1-g^{(\alpha, \beta)}_{v,v'}(a)\partial_a (g^{(\alpha, \beta)}_{v',v}(a)\beta_{v',v}(a))\alpha_{v',v}(a)\Big).\nonumber
\end{eqnarray}
Now, substituting that into the equation (\ref{aprime}) for $a'$ we obtain:
\begin{eqnarray}
a'+\partial_a f^{(\alpha, \beta)}_{v',v}(a)\Big(s_{v,v'}(a)\zeta-\alpha_{v',v}(a)\Big)\beta_{v',v}(a)=f^{(\alpha, \beta)}_{v'v}(a)\nonumber,
\end{eqnarray}
which is equivalent to 
\begin{eqnarray}
&&a'=f^{(\alpha, \beta)}_{v'v}(a)-\nonumber\\
&&\Big(\partial_af_{v'v}(a)\Big(1-\partial_a\beta_{v',v}(a)\alpha_{v',v}(a))\Big)\zeta- \partial_a f^{(\alpha, \beta)}_{v',v}(a)\alpha_{v',v}(a)\Big)\beta_{v',v}(a),\nonumber
\end{eqnarray}
and simpler:
\begin{eqnarray}
&&a'=f^{(\alpha, \beta)}_{v'v}\Big(a-(1-\partial_a\beta_{v',v}(a)\alpha_{v',v}(a))\zeta\beta_{v',v}(a)+\beta_{v',v}(a)\alpha_{v',v}(a)\nonumber
\Big).
\end{eqnarray}
One can see from the transformations we obtained that the self-dual $(1|1)$ supermanifolds are indeed $N=1$ SRS. 
Let us combine all that in the following theorem.
\begin{Thm}
Given the coordinate transformations (\ref{gnsplit2}) for $SF$, the coordinate transformations for the dual manifold $\widetilde{SF}$ of $(0|1)$ divisors is given by the formulas:
\begin{eqnarray}
&&\zeta'=g^{(\alpha, \beta)}_{v,v'}(a)\partial_a f^{(\alpha, \beta)}_{v'v}(a)\Big(s_{v,v'}(a)\zeta- \alpha_{v',v}(a)\Big), {~\rm where}\\
&&s_{v,v'}=
\Big(1-g^{(\alpha, \beta)}_{v,v'}(a)\partial_a (g^{(\alpha, \beta)}_{v',v}(a)\beta_{v',v}(a))\alpha_{v',v}(a)\Big),\nonumber\\
&&a'=f^{(\alpha, \beta)}_{v'v}\Big(a-(1-\partial_a\beta_{v',v}(a)\alpha_{v',v}(a))\zeta\beta_{v',v}(a)+\beta_{v',v}(a)\alpha_{v',v}(a)\nonumber
\Big).
\end{eqnarray}
\end{Thm}
\noindent {\bf Remark.} Note, that in the case of a dual manifold, 
$\mathcal{L}$ is replaced by $\mathcal{L}^{-1}\otimes T$. 

\section{N=2 Super Riemann Surfaces}

In this section we write down the coordinate transformations
for punctured $N=2$ supermanifold $SF_{N=2}$, corresponding to $SF$, based on the equivalence between complex structures on $(1|1)$-supermanifolds and  superconformal structures on $N=2$ supermanifolds discovered in \cite{schwarz}.  

Let us write down the transition functions between the chart with coordinates $(z,\theta)$ and chart with coordinates $(u,\eta)$ on $(1|1)$ supermanifold in the following way: 
\begin{eqnarray}
u=S(z)+\theta V(z)\varphi (z), \quad \eta=\psi(z)+\theta V(z),
\end{eqnarray}
where $S(z), V(z)$ and $\varphi (z), \psi(z)$ are correspondingly even and odd analytic functions.
Onthe other hand, the superconformal coordinate transformations for $N=2$ SRS between the charts with coordinates  $(z, \theta_{+}, \theta_{-})$ 
and $(z', \theta'_{+}, \theta'_{-})$ are:
\begin{eqnarray}
&&z'=q(z)+\frac{1}{2}\theta_-\epsilon_+(z)q_-(z)+\frac{1}{2}\theta_+\epsilon_-(z)q_+(z)+\frac{1}{4}\theta_{+}\theta_-\partial_z(\epsilon_+(z)\epsilon_-(z))\nonumber \\
&&{\theta}'_{+}=\epsilon_{+}(z)+\theta_{+}q_{+}(z)+\frac{1}{2}\theta_+\theta_-\partial_z\epsilon_{+}(z)\nonumber\\
&&{\theta}'_{-}=\epsilon_{-}(z)+\theta_{-}q_{-}(z)+\frac{1}{2}\theta_-\theta_+\partial_z\epsilon_{-}(z)\nonumber\\
&&q_{+}(z)q_{-}(z)=\partial_z q(z)+\frac{1}{2}(\epsilon_{+}(z)\partial_z\epsilon_{-}(z)+\epsilon_{-}(z)\partial_z\epsilon_{+}(z)).
\end{eqnarray} 
There is a following Theorem matching these transformations.
\begin{Thm}\cite{schwarz}
There is one-to one correspondence between $N=2$ SRS from $(1|1)$-supermanifolds. The explicit correspondence between transition functions is given by the following formulas:
\begin{eqnarray}
&&\epsilon_{+}(z)=\psi(z), \quad q_{+}(z)=V(z)\nonumber\\
&&\epsilon_{-}(z)=\varphi(z), \quad q_{-}(z)=(\partial_z S(z)-\partial_z\psi(z)\varphi(z)) V^{-1}(z),\nonumber\\
&&q(z)=S(z)+\frac{1}{2}\varphi(z)\psi(z).
\end{eqnarray} 
\end{Thm}
Let us now describe how it works for the transition functions we introduced in the previous section. 
In our case:
\begin{eqnarray}
&&V(w)=g^{(\alpha, \beta)}_{v',v}(w),\quad  \psi(w)=g^{(\alpha, \beta)}_{v',v}(w)\beta_{v',v}(w),\\
&&S(w)= f^{(\alpha, \beta)}_{v',v}(w) , \quad        
\varphi(w)=\partial_{w}f^{(\alpha, \beta)}_{v,'v}(w)\alpha_{v',v}(w)g^{(\alpha, \beta)}_{v,v'}(w). \nonumber
\end{eqnarray}
Therefore, we can write for the transition functions of $SF_{N=2}$:
\begin{eqnarray}
&&\epsilon_{+}(w)=g^{(\alpha, \beta)}_{v',v}(w)\beta_{v',v}(w),\nonumber\\
&&\epsilon_{-}(w)=\partial_{w}f^{(\alpha, \beta)}_{v',v}(w)\alpha_{v',v}(w)g^{(\alpha, \beta)}_{v,v'}(w), \nonumber\\
&& q_{+}(w)=g^{(\alpha, \beta)}_{v'v}(w),\nonumber\\
&&q_{-}(w)=\Big(\partial_w f^{(\alpha, \beta)}_{v'v}(w)-\nonumber\\
&&\partial_w(g^{(\alpha, \beta)}_{v',v}(w)\beta_{v',v}(w))\partial_{w}f^{(\alpha, \beta)}_{v',v}(w)\alpha_{v',v}(w)g^{(\alpha, \beta)}_{v,v'}(w)\Big) g^{(\alpha, \beta)}_{v,v'}(w),\nonumber\\
&&q(w)=f^{(\alpha, \beta)}_{v',v}(w) +\frac{1}{2}\partial_{w}f^{(\alpha, \beta)}_{v,'v}(w)\alpha_{v',v}(w)\beta_{v',v}(w).\nonumber
\end{eqnarray} 
This can be rewritten in a simpler way:
\begin{eqnarray}\label{n2gluing}
&&\epsilon_{+}(w)=g^{(\alpha, \beta)}_{v',v}(w)\beta_{v',v}(w),\nonumber\\
&&\epsilon_{-}(w)=\partial_{w}f^{(\alpha, \beta)}_{v',v}(w)\alpha_{v',v}(w)g^{(\alpha, \beta)}_{v,v'}(w), \nonumber\\
&& q_{+}(w)=g^{(\alpha, \beta)}_{v',v}(w)\nonumber\\
&&q_{-}(w)=\partial_w f^{(\alpha, \beta)}_{v',v}(w)g^{(\alpha, \beta)}_{v,v'}(w)(1+\alpha_{v',v}(w)\partial_w\beta_{v',v}(w))+\nonumber\\
&&\partial_wg^{(\alpha, \beta)}_{v,v'}(w)\partial_{w}f^{(\alpha, \beta)}_{v',v}(w)\beta_{v',v}(w)\alpha_{v',v}(w),\nonumber\\
&&q(w)=f^{(\alpha, \beta)}_{v',v}(w)+\frac{1}{2}\partial_{w}f^{(\alpha, \beta)}_{v,'v}(w)\alpha_{v',v}(w)\beta_{v',v}(w).
\end{eqnarray} 
Hence we obtain the following theorem.
\begin{Thm}  Formulas (\ref{n2gluing}) produce the transition functions describing the superconformal structure on $N=2$ SRS with punctures, corresponding to $(1|1)$-supermanifolds with transition functions (\ref{gnsplit2}). Namely the transition function corresponding to oriented edge $v,v'$ of the fatgraph, i.e. the overlap $U_v\cap U_{v'}$ is desribed by the functions $\epsilon_{\pm}(w)$, $q_{\pm}(w)$ from (\ref{n2gluing}).
\end{Thm}
\section{Involution and $N=1$ Super-Riemann surfaces with NS and R punctures.}
\subsection{Involution: R vs NS punctures}
There is an involution $I$ on the moduli space of super-Riemann surfaces, such that 
\begin{eqnarray}
I: D_{\pm}\to D'_{\mp},
\end{eqnarray}
where $D'_{\mp}$ is the corresponding operator after $N=2$ superconformal transformation.

Such an involution takes $N=2$ super Riemann surface to the {\it dual}, which on the level of $(1|1)$- supermanifolds produces a manifold of $(0|1)$-divisors, which we discussed earlier.  The self-dual supermanifolds are known to be $N=1$ super-Riemann surfaces. 

Let us describe how this works on a $N=2$ supertube (or $N=2$ punctured disk) with coordinates $( x, \eta_+, \eta_-)$, where 
$x\sim x+2\pi i$. 
Let us consider an obvious choice of how involution could act in these coordinates:
\begin{eqnarray}
D_{+}\to D_{ -}, ~ D_{-}\to D_{+}, 
\end{eqnarray} 
For self-duality one has to 
identify $\eta_{+}$ and $\eta_-$, i.e. $(x, \eta_{+}, \eta_{-})\sim (x +2\pi, \eta_{-}, \eta_{+})$. The operator $D=D_++D_-$ gives a standard superconformal structure on a supertube. We see, that in this case the puncture is a {\it Ramond puncture}.
Let us perform an elementary $N=2$ superconformal transformation, amounting to reflection, so that involution is 
\begin{eqnarray}
D_{+}\to -D_{ -}, ~ D_{-}\to -D_{+}, 
\end{eqnarray} 
i.e. $\eta_{\pm}\to -\eta_{\mp}$. The invariance under this involution gives the identification $(x, \eta_{+}, \eta_{-})\sim (x +2\pi, -\eta_{-}, -\eta_{+})$, so that  operator $D=D_++D_-$ gives a superconformal structure around NS puncture. 

Note, that the two examples of the action of involution which we considered in this section are the only ones, which preserve the base manifold. 

\subsection{Split $N=1$ SRS}
Let us now discuss split $N=2$ SRS, which implies that we let cocycles $\alpha,\beta=0$. 
The involution
\begin{eqnarray}
I: D_{\pm}\to D_{\mp},
\end{eqnarray}
acts on the level of transition functions as follows:
$$
q_{\pm}(z)\to q_{\mp}(z)
$$
Therefore, for fixed points of the involution we have 
\begin{eqnarray}
g_{v',v}^2(w)=\partial_w f_{v',v}(w).
\end{eqnarray}
This means that $g^{(\alpha, \beta)}_{v',v}(w)={\rm sign}(v',v) \sqrt{\partial _w f_{v',v}(w)}$ where ${\rm sign}(v',v)$ is the notation for the sign of the square root, so that on a resulting $N=1$ SRS 
we have: 
\begin{eqnarray}
\xi'={\rm sign}(v',v)\sqrt{\partial_w f_{v',v}(w)}\xi.
\end{eqnarray}
Choice of signs for such square roots is the same as the choice of spin structure on the punctured surface. However, we already discussed that problem on the level of fatgraphs (Section 3), which leads to the following Theorem.

\begin{Thm}\label{splitsrs}
Consider a metric fatgraph $\tau$ with a spin structure $\omega$ provided by the orientation as discussed in Section 3.  This data defines the superconformal structure on the split $N=1$ SRS. 
For every boundary cycle on the fatgraph, corresponding to puncture $p$, let $m_p$ be the number of oriented edges, which are opposite to the orientation induced by the one on the surface. The corresponding puncture is Ramond or Neveu-Schwarz, depending on whether $m_p$ is even or odd. 
\end{Thm}
\noindent {\bf Proof.}
So, let us consider the metric graph $\tau$ with orientations on edges. Our problem is to use orientations to define 
To do that, for each overlap we will look at the $z$ coordinates on stripes, discussed in section 4. For given vertices $v$ and $v'$, the transformation between $z$ and $z'$ coordinates is given by 
$$
z'=\tilde{f}_{v',v}(z)=L_{v,v'}-z
$$
We will define the value of the $\sqrt{\partial_z \tilde{f}_{v',v}(z)}=\pm i$ in the following way. If the orientation is from vertex $v$ to $v'$ we choose the positive sign (${\rm sign}(v,v')=1$) and negative otherwise (${\rm sign}(v,v')=-1$). One can  prove that such choice does not depend on the choice of orientation for a given spin structure, namely a different choice, corresponding to a fatgraph reflection, will just result in a reflection of an odd coordinate for a given vertex. 

Regarding R and NS punctures, one can deduce immediately that the statement is correct just by a simple condition that there is a natural combinatorial constraint on the punctures with $m_p$ being odd on a fatgraph (see Section 4), matching the one for Ramond punctures on a surface. 
Nevertheless, let us prove that directly. 

For a given choice of spin structure, let us superconformally continue $g^{(\alpha, \beta)}_{v'v}$ cocycles by constructing $\tilde {g}_{p,v}(z)={\rm sign}(p,v)\sqrt{\partial_z\tilde{f}_{p,v}(z)}$ on 
$U_v\cap U_{v'}\cap V_p$, 
where we remind that  $$x=\tilde{f}_{p,v}(z)=\frac{2\pi i}{a_B}(L_{v_1,v_2}+...+ L_{v_k,v_{k-1}}+z),$$ 
where $z$ is the coordinate on the consequtive  stripe 
$v_k,v_{k+1}$ and $v_1,\dots , v_m$ are consequtive vertices around the puncture.
Now we obtain $R$ and $NS$ punctures by gluing the supertube with a proper twist of the odd variable. That will of course depend on the number $m_p$ of the edges $\{v_i, v_{i+1}\}$, which  have opposite orientation with respect to orientation induced on the cycle by the one on the surface.
Note, that in terms of $z$-variables $\sqrt{\partial_z\tilde{f}_{p,v}(z)}$ is a constant, so one can again make a choice of signs explicitly. We have the following identity 
\begin{eqnarray}\label{cyc}
{\rm sign}(v_1,v_2){\rm sign}(v_2,v_3)\dots {\rm sign}(v_{n-1},v_n){\rm sign}(v_n,v_1)=\pm 1,
\end{eqnarray}
where positive sign is for even $N$ and negative otherwise. In the case of $m_p$ even, we can choose $\{{\rm sign}(p,v)\}$ so that for ${\rm sign}(v,v')$, so that $v,v'$ are neighboring vertices, so that  ${\rm sign}(v,v')={\rm sign}(p,v){\rm sign}(p,v')$, thus gluing the stripes into supertube. 
However, this is not possible in the case 
of odd $m_p$. In this case we have to assume that 
${\rm sign}(v_n,v_1)=-{\rm sign}(p,v_1){\rm sign}(p,v_n)$, thus gluing the stripes into twisted 
supertube corresponding to $NS$ puncture. \hfill $\blacksquare$\\
% tranformations between fermion variables on $U_v\cap V_p$ is: 
 %$$\xi={\rm sign}(v,p)\sqrt{\partial _x f_{v,p}(x)}\eta.$$ 

\noindent {\bf Remark.} One can of course superconformally transform the twisted supertube in $NS$ puncture case into the disk, the same way we did in the introduction, thus making 
the corresponding cocycle $\{g_{v,p}(y)\}=\{\pm\sqrt {\partial _y f_{v,p}(y)}\}$ . In Ramond case this is of course impossible. We see that if $p$ is an R puncture, $$g_{v,p}^2(x)=y\partial _y f_{v,p}(y).$$ Therefore, for the bundle $\mathcal{L}$ we have a condition:
\begin{eqnarray}\label{bundle}
\mathcal{L}^2=T\otimes \mathcal{O}\Big(-\sum^{n_R}_{i=1}{p_i}\Big), 
\end{eqnarray}
which is possible only when $n_R$  is divisible by 2.

\subsection{$N=1$ SRS: non-split case.}
%In order to find the explicit formulas for the global transition functions non-split $N=1$ SRS, one has to relax certain properties for the formulas in (\ref{n2gluing}) in order to impose the invariance under simple involution $D_{\pm}\to D_{\mp}$. In other words, we have to assume $g^{(\alpha, \beta)}_{v',v}$, $\alpha_{v',v}$
%$\beta_{v',v}$ absolutely random analytic functions on the overlaps $U_v\cap U_v'$. 

%\begin{equation}
%g^{(\alpha, \beta)}_{v',v}=\pm \sqrt{\partial_{w}f_{v',v}(w)}(1+\frac{1}{2}\alpha_{v',v}(w)\partial_w\alpha_{v',v}(w)), 
%\end{equation}

In order to construct nonsplit $N=1$ SRS, we first will do it on infinitesimal level near the split $N=1$ SRS. So, let us look at the formulas (\ref{n2gluing}) when $\alpha_{v',v}$
$\beta_{v',v}$ are infinitesimal: 
\begin{eqnarray}
&&\epsilon_{+}(w)=g_{v',v}(w)\beta_{v',v}(w),\nonumber\\
&&\epsilon_{-}(w)=\partial_{w}f_{v',v}(w)\alpha_{v',v}(w)g^{(\alpha, \beta)}_{v,v'}(w), \nonumber\\
&& q_{+}(w)=g_{v',v}(w),\nonumber\\
&&q_{-}(w)=\partial_w f_{v',v}(w)g_{v,v'}(w),\nonumber\\
&&q(w)=f_{v',v}(w).
\end{eqnarray} 
An invariance under simple involution $D_{\pm}\to D_{\mp}$ allows to identify  $\alpha_{v',v}$ and 
$\beta_{v',v}$ and as before ${g^2_{v',v}}(w)=\partial_w f_{v',v}(w)$, thus infinitesimally the transformations for the resulting $N=1$ SRS  on the overlap $U_v\cap U_{v'}$ is given by:
\begin{eqnarray}\label{tesm}
&&w'=f_{v'v}(w+\xi\rho_{v',v}(w))\nonumber\\
&&\xi'=\pm\sqrt{\partial_wf_{v',v}(w)}(\xi+\rho_{v',v}(w)),
\end{eqnarray} 
so that the signs of are prescribed as in the  Theorem \ref{splitsrs}, 
where  $$\rho \in \Pi\check{Z}^{1}(F^c, \mathcal{L}\otimes O(-D_{NS}))~ {\rm and} ~\mathcal{L}^2=T\otimes O(-D_R),$$ so that $D_R$ and $D_{NS}$ are the divisors corresponding to the sum of all $NS$ and $R$ punctures correspondingly. We described such cocycles using odd number decorations at the vertices of the  the fatgraph in the Theorem \ref{cycdes}.   
%\begin{equation}
%g^{(\alpha, \beta)}_{v',v}=\pm \sqrt{\partial_{w}f_{v',v}(w)}(1+\frac{1}{2}\alpha_{v',v}(w)\partial_w\alpha_{v',v}(w)), 
%\end{equation}
The formulas (\ref{tesm}) are not hard to continue to full superconformal transformations for transition functions (one can obtain them by applying involution $D_{\pm}\to D_{\mp}$ invariance to the formulas (\ref{n2gluing}) as well):
\begin{eqnarray}
&&w'=f^{(\rho)}_{v'v}(w+\xi\lambda^{(\rho)}_{v',v}(w))\nonumber\\
&&\xi'=\pm\sqrt{\partial_w f^{(\rho)}_{v',v}(w)}(1+\frac{1}{2}\lambda^{(\rho)}_{v',v}(w)\partial_w\lambda^{(\rho)}_{v',v}(w))(\xi+\lambda_{v',v}(w)).\nonumber
\end{eqnarray}

Combining the parametrization data for cocyles $\rho$ from Theorem \ref{cycles} with the results of this section, we obtain the following omnibus Theorem, describing the dense set of superconformal structures oinside moduli space of $N=1$ SRS.

\begin{Thm}\label{maintheorem}

\item Consider the following data on a fatgraph $\tau$:
\begin{enumerate}
\item Metric structure.\\

\item Spin structure, as equivalence class of orientations on the fatgraph. The cycles on the fatgraph, encircling the punctures are divided into two subsets, $NS$ and $R$, depending on whether there is odd or even number of edges oriented opposite to the surface-induced orientation of the appropriate boundary piece of a fatgraph correspondingly. We denote the number of the corresponding boundary pices as $n_R$ and $n_{NS}$. \\

\item Ordered  set $\{\sigma_v^k\}_{k=0,\dots, m_v-3}$ of odd complex parameters for each vertex $v$, where $m_v$ is the valence of the vertex $v$. \\

Then the following is true:\\

\begin{enumerate}
\item Data from {\rm (1)} and {\rm (2)} determine uniquely the split Riemann surface with $n_R$ Ramond and $n_{NS}$ Neveu-Schwarz punctures with the transition functions given by 
\begin{eqnarray}
w'=f_{v',v}(w)\quad \xi'=\pm\sqrt{\partial_w f_{v',v}(w)}\xi,
\end{eqnarray}
one for each overlap $U_v\cap U_{v'}$.  The sign of the square root is given by the spin structure on the fatgraph, making odd coordinate a section of a line bundle $\mathcal{L}$ on the corresponding closed Riemann surface $F^c$, such that $\mathcal{L}^2=T\otimes \mathcal{O}(-D_R)$, where $D_R$ is a divisor, which is a sum of points corresponding to the Ramond punctures.\\

\item Part {\rm (3)}  of the above data  allows to construct \v{C}ech cocycles on a  Riemann surface F, which  are the representatives of $\Pi\check{H}^1(F^c, \mathcal{L}\otimes \mathcal{O}(-D_{NS}))$, where $D_{NS}$ is a divisor, corresponding to the sum of the points corresponding to NS punctures: 
\begin{eqnarray}
&&\rho_{v,v'}|_{U_v\cap U_{v'}}=\rho_v-\rho_{v'},~{\rm  so~ that ~ } ~\rho_v|_{U_v}=\frac{\sigma_v(w)}{w^{m_v-2}},~  \rho_{v'}|_{U_v'}=\frac{\sigma_{v'}(w')}{w'^{~m_{v'}-2}},\nonumber \\
&&\sigma_v(w)=\sum^{m_{v}-3}_{i=0}\sigma^i_v w^i, \quad  \sigma_{v'}(w')=\sum^{m_{v'}-3}_{i=0}\sigma^i_{v'} {{w'}}^i,
\end{eqnarray}
where $\rho_v$, $\rho_{v'}$ are the meromorphic sections of 
$\mathcal{L}\otimes \mathcal{O}(-D_{NS})$ on $U_v$, $U_{v'}$ correspondingly,  so that $m_v$ is valence of the given vertex $v$.

The cocycles defined by configurations described by $\{\sigma_v\}$  and $\{\tilde{\sigma}_v\}$ are equivalent to each other if and only if 
\begin{eqnarray}\label{equival}
\sigma_v(w)-\tilde{\sigma}_v(w)=\gamma^{(m-3)}(w), 
\end{eqnarray}
for every vertex $v$,  
$\gamma\in \Pi{H}^0(F^c, \mathcal{L}\otimes K^2\otimes \mathcal{O}(D_{NS}+2D_{R}))$,   
$\gamma|_{U_v}=\gamma(w)$ so that $\gamma^{(m-3)}(w)$ is the Taylor expansion of $\gamma(w)$ up to order $m-3$.\\

We call two sets of data associated to the fatgraph $\tau$ equivalent, if they are related as in (\ref{equival}).\\

\item There exist a  superconformal structure  for $N=1$ super Riemann surface SF with $n_R$ Ramond punctures and $n_{NS}$ Neveu-Schwarz punctures so that the superconformal transition functions on for each overlap $U_v\cap U_{v'}$  are:
\begin{eqnarray}
&&w'=f^{(\sigma)}_{v'v}(w+\xi \lambda^{(\sigma)}_{v',v}(w))\nonumber\\
&&\xi'=\pm\sqrt{\partial_w f^{(\sigma)}_{v',v}(w)}\Big(1+\frac{1}{2}\lambda^{(\sigma)}_{v',v}(w)\partial_w\lambda^{(\sigma)}_{v',v}(w)\Big)(\xi+\lambda^{(\sigma)}_{v',v}(w)),
\end{eqnarray}
where the deformed functions $f^{(\sigma)}_{v',v}$, $\lambda^{(\sigma)}_{v',v}$ depend on odd parameters $\{\sigma^{k}_v\}$, characterizing \v{C}ech cocylce $\{\rho_{v',v}\}$,  with $f^{(0)}_{v',v}=f_{v',v}$ and in the first order in $\{\sigma^{k}_v\}$ variables $\lambda^{(\sigma)}_{v',v}=\rho_{v',v}$. 
 \\
 
\item 
To describe the non-split SRS, we fix the choice of transition functions in ${\rm (c)}$ for every metric spin fatgraph $\tau$ with the odd data from  {\rm (3)}. We consider the set of superconformal structures constructed by picking one superconformal structure per equivalence class of data for every fatgraph $\tau$. The points in this set represent inequivalent superconformal structures, and together they form a dense subspace of odd complex dimension $2g-2+n_{NS}+n_{R}/2$  in the space of all superconformal structures with $n_{NS}$ Neveu-Schwarz and $n_{R}$ Ramond punctures associated to $F$.  

\end{enumerate}
\end{enumerate}
%In other words, Strebel differentials with the prescribed odd data at their zeroes and $NS$ or $R$ assignment at poles modulo the described above equivalences from ii) parametrize the $\mathbb{R}^s_{+}$-bundle over dense set in the corresponding supermoduli space of $N=1$ SRS.
\end{Thm}

\end{document}